\documentclass [12pt]{article}
\title {On the forking topology of a reduct of a simple theory}
\author {Ziv Shami\\Ariel University}
\newtheorem {theorem}{Theorem}[section]
\newtheorem {lemma}[theorem]{Lemma}
\newtheorem {definition}[theorem]{Definition}

\newtheorem {fact}[theorem]{Fact}

\newtheorem {corollary}[theorem]{Corollary}

\newtheorem {remark}[theorem]{Remark}

\newtheorem {proposition}[theorem]{Proposition}

\newtheorem {claim}[theorem]{Claim}

\newtheorem {example}[theorem]{Example}
\newtheorem {example-revisited}[theorem]{Example-revisited}

\def\proof {\noindent \textbf{Proof:} }

\def\qed {$\ \ \ \ \Box$}

\newsavebox{\indbin}
\savebox{\indbin}{\begin{picture}(0,0)
\newlength{\gnu}
\settowidth{\gnu}{$\smile$} \setlength{\unitlength}{.5\gnu} \put(-1,-.65){$\smile$}
\put(-.25,.1){$|$}
\end{picture}}
\newcommand{\nonfork}[3]
{\mbox{$\begin{array}{ccc} \mbox{$#1$} & \usebox{\indbin} & \mbox{$#2$} \\
        & \mbox{$#3$} &
\end{array}$}}
\newcommand{\nonforkempty}[2]
{\mbox{$\begin{array}{ccc} \mbox{$#1$} & \usebox{\indbin} & \mbox{$#2$}
\end{array}$}}

\newsavebox{\sindbin}
\savebox{\sindbin}{\begin{picture}(0,0)
\newlength{\sgnu}
\settowidth{\sgnu}{$\smile$} \setlength{\unitlength}{.5\sgnu} \put(-1,-.65){$\smile$}
\put(-.25,.1){$|s$}
\end{picture}}

\newsavebox{\minusindbin}
\savebox{\minusindbin}{\begin{picture}(0,0)
\newlength{\minusgnu}
\settowidth{\minusgnu}{$\smile$} \setlength{\unitlength}{.5\minusgnu} \put(-1,-.65){$\smile$}
\put(-.25,.1){$|{^-}$}
\end{picture}}
\newcommand{\minusnonfork}[3]
{\mbox{$\begin{array}{ccc} \mbox{$#1$} & \usebox{\minusindbin} & \mbox{$#2$} \\
        & \mbox{$#3$} &
\end{array}$}}
\newcommand{\minusnonforkempty}[2]
{\mbox{$\begin{array}{ccc} \mbox{$#1$} & \usebox{\minusindbin} & \mbox{$#2$}
\end{array}$}}

\newsavebox{\qindbin}
\savebox{\qindbin}{\begin{picture}(0,0)
\newlength{\qgnu}
\settowidth{\qgnu}{$\smile$} \setlength{\unitlength}{.5\qgnu} \put(-1,-.65){$\smile$}
\put(-.25,.1){$|_{qf}$}
\end{picture}}

\def\card #1 {{\vert #1 \vert}}

\def\CC {{\cal C}}

\def\UU {{\cal U}}

\def\MM {{\cal M}}

\usepackage{amsfonts}

\begin{document}
\maketitle

\begin{abstract}
Let  $T$ be a simple $L$-theory and let $T^-$ be a reduct of $T$ to a sublanguage $L^-$ of $L$. For variables $x$, we call an $\emptyset$-invariant set $\Gamma(x)$ in $\CC$ a \em universal transducer \em if  for every formula $\phi^-(x,y)\in L^-$ and every $a$,  $$\phi^-(x,a)\ L^-\mbox{-forks\ over}\ \emptyset\ \mbox{\ iff\ } \Gamma(x)\wedge \phi^-(x,a)\ L\mbox{-forks\ over}\ \emptyset.$$ We show that there is a greatest universal transducer $\tilde\Gamma_x$ (for any $x$) and it is type-definable. In particular, the forking topology on $S_y(T)$ refines the forking topology on $S_y(T^-)$ for all $y$. Moreover, we describe the set of universal transducers in terms of certain topology on the Stone space and show that $\tilde\Gamma_x$ is the unique universal transducer that is $L^-$-type-definable with parameters. If $T^-$ is a theory with the wnfcp (the weak nfcp) and $T$ is the theory of its lovely pairs of models we show that $\tilde\Gamma_x=(x=x)$ and give a more precise description of the set of universal transducers for the special case where $T^-$ has the nfcp.
\end{abstract}

\section{Introduction}
The forking topology for simple theories, introduced in [S], is a generalization of topologies introduced by Hrushovski [H0] and Pillay [P].  It is the minimal topology on $S_x(A)$ such that all the
relations $\Gamma_F(x)$ defined by $\Gamma_F(x)=\exists y (F(x,y)\wedge \nonfork{y}{x}{A})$ are
closed for any type-definable relation $F(x,y)$ over $A$.

Originally, a version of this topology has been introduced (around 1984) by Hrushovski  [H0] for the (unpublished) proof of superstability of countable unidimensional stable theories; in the proof, an unbounded set of finite rank is constructed that is open in the forking topology.
In [P], where supersimplicity of any countable unidimensional wnfcp hypersimple theory (i.e. a simple theory that eliminates hyperimaginaries) is established, the topology has been modified to work for theories with the weak non finite cover property (wnfcp), an analogue of the non finite cover property (nfcp) for simple theories.  In [S] we modified the topology defined in [P] and proved general theorems for simple theories related to unidimensional theories. The forking topology turned out to be quite a powerful tool and had several applications: finite length analysis of any type analyzable in a forking open set provided that the forking topologies are closed under projections (e.g. $T$ has wnfcp) [S],  supersimplicity of countable (and large class of uncountable)  unidimensional  hypersimple theories [S1,S2] and a generalization of Buechler's dichotomy for $D$-rank 1 types to simple theories [S3].

In this paper, we fix a simple $L$-theory $T$ and a reduct $T^-$ of $T$ to a sublanguage $L^-$ and present a way in which the forking topology of $T^-$ can be recovered from the forking topology of $T$; it is done through the notion of a universal transducer that is defined in the abstract in a restrictive form (a more general setting is presented in the paper).
Moreover, we characterize the set of universal transducers via a new topology, we call the $NFI$-topology (see Definition \ref{def. NFI-top}). Our observation that the forking-topology of a simple theory refines the forking-topology of any reduct is, in way, a substitute for the fact that forking-independence in a simple theory does not, in general, strengthen forking-independence in a reduct. Our hope is to find more relationships between forking-independence in a simple theory and forking-independence in a reduct; moreover, we expect to find more connections between sets related to the forking topology (e.g sets defined by the NFI-topology) and sets that are both $L$-type-definable over $\emptyset$ and $L^-$-type definable with parameters. In particular, we expect that the $NFI$-topology could be proved to be  $L^-$-invariant over parameters for any simple theory (Lemma \ref{lemma2}(2) confirms this for stable theories).  The following is a special case of our main theorems.
\begin{theorem}
Given variables $x$, there is a greatest (with respect to inclusion) $\emptyset$-invariant subset of $\CC^x$ that is a universal  transducer. Denote this subset by $\tilde\Gamma_x$. Then, $\tilde\Gamma_x$ is $L$- type-definable and it is the unique universal transducer that is $L^-$-type definable with parameters. If $T$ is stable, we have the following characterization of the set of universal transducers: an $\emptyset$-invariant set  $\Gamma(x)$ in $\CC$  is a universal transducer iff $\Gamma(x)$ is a dense subset of  $\tilde\Gamma_x$ in the relative topology on $\tilde\Gamma_x$ generated by the family of $L$-formulas $\phi(x)$ over $\emptyset$ that are $L^-$-definable with parameters.
\end{theorem}

In particular, the reduct map from $S_y(T)$ to $S_y(T^-)$ (for any variables $y$) is continuous with respect to the forking topologies on the Stone spaces.  Lastly, we get a more precise information in the special case of lovely-pairs: we look at the case where the reduct theory ($T^-$ in our general setting)  is an arbitrary theory with the wnfcp,  denoted by $T$ (in a language $L$), and at the expansion of it  $T_P$ ($T$ in our general setting) defined as the theory of its lovely pairs of models (in the language $L_P=L\cup \{P\}$). The result we obtained for $T_P$ and the reduct $T$ is the following.

\begin{proposition}
For any variables $x$, $\tilde\Gamma_{x}=(x=x)$, namely the greatest universal transducer in the variables $x$ is  $(x=x)$.
If $T$ is in addition stable (equivalently $T$ has nfcp), then an $L_P$--invariant set over $\emptyset$
is a universal transducer iff it intersect every non-empty $L$-definable set over $\emptyset$.
\end{proposition}

We assume basic knowledge of simple theories as in [K],[KP],[HKP]. A good textbook on simple theories is [W]. In this paper, unless otherwise stated,  $T$ will denote a complete first-order simple theory in an arbitrary language $L$ and we work in a $\lambda$-big model $\CC$ of $T$ (i.e. a model with the property that any expansion of it by less than $\lambda$ constants is splendid) for some large $\lambda$. We call $\CC$ the monster model. Note that any $\lambda$-big model (of any theory) is $\lambda$-saturated and $\lambda$-strongly homogeneous and that $\lambda$-bigness is preserved under reducts (by Robinson consistency theorem). We use standard notations. For a small subset $A\subseteq\CC$, $T_A$ will denote the theory of $(\CC,A)$ ($\CC$ expanded by constants for each $a\in A$). Partial types are usually identified with the set of their solutions in the monster model. For an invariant set of a fixed sort (or finitely many) we write (e.g.) $\UU(x)$ where $x$ is a finite tuple of variables suitable for these sorts. For variables $x$, $\CC^x$ denotes the set of tuples from $\CC$ whose sort is the sort of $x$. An invariant set of possibly some distinct sorts will be denoted by (e.g.) $\UU$ (with no variables added). If $\UU$ is a set we denote by $\UU^{<\omega}$ the set of all finite sequences of elements in $\UU$. For a partial type $p$ over a model, $Cl(p)$ denotes the set of formulas $\phi(x,y)\in L$ that are represented in $p$.

\section{Transducers}
In this section we introduce the notion of a universal $F$-transducer for an $\emptyset$-invariant set $F$ and prove generalizations of the results stated in the abstract for a simple theory and a reduct. First, recall the definition of the forking topology.

\begin{definition}\em [S, Definition 2.1]\em  \label{tau definition}
\em Let $A\subseteq\CC$ and let $x$ be a finite tuple of variables. A set $U\subseteq S_x(A)$ is said to be \em a basic forking-open set over $A$ \em if
 there exists $\phi(x,y)\in L(A)$ such that $$U=\{p\in S_x(A) \vert\ \phi(a,y)\ \mbox{forks\ over}\ A\ \mbox{for\ all\ } a\models p \}.$$
\end{definition}

\noindent We identify subsets of $S_x(A)$ with $A$-invariant sets. Note that the family of basic forking-open sets over $A$ is closed under finite
intersections, thus form a basis for a unique topology on $S_x(A)$ which we call the
forking-topology or the forking-topology.

\begin{remark} \label{ft_remark}\em
Note that the forking-topology on $S_x(A)$ refines the Stone-topology (for every $x$ and $A$) and
that $$\{a\in\CC^x\vert a\not\in acl(A)\}(=\{a\in\CC^x\vert x=a\ \mbox{forks\ over}\ A\})$$ is a
forking-open subset of $S_x(A)$.
\end{remark}

We fix now the notations for the rest of this section. $T^-$ will denote a reduct of $T$ to some sublanguage $L^-$ of $L$, i.e. $T^-$ is the set of $L^-$-sentences in $T$. We will assume for simplicity of notation that $L^-$ and $L$ have the same set of sorts (the general case is very similar and discussed in Remark \ref{neq_sorts}). Let $\CC^-=\CC\vert L^-$. As mentioned in the introduction, we know that both $\CC$ and $\CC^-$ are highly saturated and highly strongly-homogeneous. $\CC^{heq}$  will denote the set of hyperimaginaries of small length ($<\lambda$) of $\CC$ and $\CC^{heq-}$ will denote the set of hyperimaginaries of small length of $\CC^-$.
We use $\nonforkempty {}{}$ to denote independence in $\CC$, and $\minusnonforkempty{}{}$ to denote independence in $\CC^-$.
For a small set $A\subseteq \CC^{heq}$, $BDD(A)$ denotes the set of countable (length)  hyperimaginaries in $\CC^{heq}$  that are in the bounded closure of $A$ in the sense of $\CC$.
For an $\emptyset$-invariant set $F$ in $\CC$ (or for a small set $F$), we denote by $bdd(F)$ ($dcl^{heq}(F)$) the set of all countable hyperimaginaries in $\CC^{heq-}$ that are in the bounded (definable) closure in the sense of $\CC^-$ of some small subset of $F$. For a small set $A\subseteq \CC^{eq}$, $ACL^{eq}(A)$ denotes the set of imaginaries in $\CC^{eq}$  that are in the algebraic closure of $A$ in the sense of $\CC$. For a small set $A\subseteq \CC^{eq-}=(\CC^-)^{eq}$, $acl^{eq}(A)$ denotes the set of imaginaries in $\CC^{eq-}$  that are in the algebraic closure of $A$ in the sense of $\CC^-$.
For a small set $X\subseteq\CC^{heq}$, let $X^-=X\cap \CC^{heq-}$.\\

Let $\Gamma(x)$ be a $B$-invariant set in $\CC$ and let $A$ be any small set. We say $\Gamma(x)$ $L$-doesn't  fork over $A$
if for some $c\models\Gamma(x)$, $\nonfork{c}{B}{A}$.\\
From now on $F$ will denote an arbitrary $\emptyset$-invariant set in $\CC$.

\begin{definition}\em
Let $\Gamma(x)$ be an $\emptyset$-invariant set in $\CC$.\\
1) We say that $\Gamma(x)$ is an \em upper universal $F$-transducer \em if for every $\bar a\in F^{<\omega}$ and $\phi^-(x,\bar y)\in L^-$, if $\Gamma(x)\wedge \phi^-(x,\bar a)$ $L$-doesn't fork over $\emptyset$, then $\phi^-(x,\bar a)$ $L^-$-doesn't fork over $\emptyset$.\\
2) We say that $\Gamma(x)$ is a \em lower universal $F$-transducer \em if for every $\bar a\in F^{<\omega}$ and  $\phi^-(x,\bar y)\in L^-$, if $\phi^-(x,\bar a)$ $L^-$-doesn't fork over $\emptyset$, then $\Gamma(x)\wedge \phi^-(x,\bar a)$ $L$-doesn't fork over $\emptyset$.\\
3)  We say that $\Gamma(x)$ is a \em universal  $F$-transducer \em if $\Gamma(x)$ is both an upper universal $F$-transducer and a lower universal $F$-transducer.\\
Whenever $F$ is omitted in 1)-3), it means $F=\CC$.
\end{definition}

\begin{example}\em
Let $T^-$ be the theory of an infinite set with no structure and let $T$ be an expansion of $T^-$ by some small set of constants $C\subseteq\CC^-$. Note that, if $x$ is a single variable,
then the type $\Gamma(x)=\{x\neq c\ \vert\ c\in C\}$ is the unique universal transducer in the variable $x$.
\end{example}

\begin{remark}\em\label{trans_fk_top}
Note that the existence of a type-definable universal transducer in any variables implies that the forking-topology of $T$ on  $S_y(T)$ refines the forking-topology of $T^-$ on $S_y(T^-)$ for every $y$, that is,
the reduct map from $S_y(T)$ to $S_y(T^-)$ (for any variables $y$) is continuous with respect to the forking topologies on these spaces: if $\Gamma(x)$ is a type-definable universal transducer over $\emptyset$ then  for every formula $\phi^-(x,y)\in L^-$, we have: $$\{b \vert\ \phi^-(x,b)\ L^-\mbox{-forks\ over}\ \emptyset\}=\bigcup_{\psi(x)\in\Gamma(x)} \{b \vert\ \psi(x)\wedge\phi^-(x,b)\ L\mbox{-forks\ over}\ \emptyset\}.$$
\end{remark}

\begin{definition}\em
For variables $x$,  we define the following $\emptyset$-invariant sets in $\CC$:\\
1) $$\tilde\Gamma_{x,F}=\{b\in \CC^x\vert\ \forall \bar a\in F^{<\omega}\ \exists \bar a'\models tp_L(\bar a) \ (\minusnonforkempty{b}{\bar a'})\}.$$
2) $$\Gamma_{x,F}^*=\{b\in \CC^x \vert\ \forall \bar a\in F^{<\omega} (\nonforkempty{b}{\bar a}\rightarrow\minusnonforkempty{b}{\bar a})\}.$$
3) $$B_{x,F}=\{b\in \CC^x \vert\ \minusnonforkempty{b}{bdd(F)\cap BDD(\emptyset)^-}\}.$$
Whenever $F$ is omitted in 1)-3), it means $F=\CC$.
\end{definition}

\begin{remark}\label {remark1}
For variables $x$, we have $$\tilde\Gamma_x=\{b\in \CC^x\vert\ \forall\phi(y)\in L:\ [\exists y\phi(y)\rightarrow \exists a\models \phi(y)\ (\minusnonforkempty{b}{a})]\}.$$
Moreover, for every model $M\models T$,
$\tilde\Gamma_x=\{b\in \CC^x\vert\ \exists M' \models tp_L(M)( \minusnonforkempty{b}{ M' } )\}$.
\end{remark}

\proof Just compactness.\qed\\

\begin{lemma}\label{lemma1}
For any variables $x$, we have $\tilde\Gamma_{x,F}=\Gamma_{x,F}^*=B_{x,F}$.
\end{lemma}

\proof  To show $\tilde\Gamma_{x,F}\subseteq B_{x,F}$ we observe:

\begin{claim}\label{claim0}
Let $M$ be a sufficiently saturated model of $T$. Then $$bdd(F)\cap BDD(\emptyset)^-=bdd(F^M)\cap BDD(\emptyset)^-.$$
\end{claim}

\proof Let $e\in bdd(F)\cap BDD(\emptyset)^-$. Then there exists  a small subset  $F_e\subseteq F$ (in fact of size at most $\vert T\vert$)
such that $e\in bdd(F_e)\cap BDD(\emptyset)^-$. Since  $M$ is sufficiently saturated, $e\in M^{heq-}$ (if $e=a/E$ then on $tp_L(a)$ there are at most $2^{{\vert T\vert }^+}$ many $E$-classes).
By saturation $M$, there exists $F'_e\subseteq M$ such that $tp_L(F'_e/e)=tp_L(F_e/e)$ and so $e\in bdd(F^M)$.\qed\\

 \noindent Now, let $b\in \tilde\Gamma_{x,F}$. By compactness, there exists a sufficiently saturated model $M'$ of $T$ such that
 $\minusnonforkempty{b}{F^{M'}}$, so  $\minusnonforkempty{b}{bdd(F^{M'})}$. By Claim \ref{claim0} we are done. To show $B_{x,F}\subseteq \Gamma_{x,F}^*$ recall the following.

\begin{fact} {\em[HN, Theorem 2.2]\em} \label{fact1}\em
Let $A,C\subseteq \CC^{heq-}$ and let $B\subseteq\CC^{heq}$ be boundedly closed in $\CC^{heq}$. Assume $\nonfork{A}{C}{B}$.
Then  $\minusnonfork{A}{C}{B^-}$.\qed
\end{fact}

\noindent Now, let $b\in B_{x,F}$ and assume $\nonforkempty{b}{\bar a}$ for some $\bar a\in F^{<\omega}$ .  By Fact \ref{fact1},  $$\minusnonfork{b}{\bar a}{BDD(\emptyset)^-}\ (*).$$  From now on work in  $\CC^-$. Let $e^-=Cb^-(Lstp(\bar a/BDD(\emptyset)^-,b))$.  $e^-$ is in the definable closure of a Morley sequence of $Lstp(\bar a/BDD(\emptyset)^-,b)$, since $\bar a\in F^{<\omega}$, we conclude $e^-\in bdd(F)$. By (*), $e^-\in BDD(\emptyset)^-$ (note that $BDD(\emptyset)^-$ boundedly closed in $\CC^{heq-}$). Thus $$\minusnonfork{_{\bar a}}{_{BDD(\emptyset)^-,b}}{_{BDD(\emptyset)^-\cap bdd(F)}}.$$ As $b\in B_{x,F}$, transitivity yields $\minusnonforkempty{b}{\bar a}$. The inclusion $\Gamma_{x,F}^*\subseteq \tilde\Gamma_{x,F}$ is immediate by extension. This completes the proof of Lemma \ref{lemma1}. \qed

\begin{proposition}\label{prop1}
For variables $x$, there exists a greatest (with respect to inclusion) $\emptyset$-invariant subset of $\CC^x$ that is a universal  $F$-transducer. Denote this subset by $\Gamma_{x,F}$.
Then, $\Gamma_{x,F}$ is also such greatest  upper universal $F$-transducer, $\Gamma_{x,F}=\tilde\Gamma_{x,F}=\Gamma^*_{x,F}$ and $\Gamma_{x,F}$ is type-definable.
In particular,  the forking-topology of $T$ on  $S_y(T)$ refines the forking-topology of $T^-$ on $S_y(T^-)$ for every $y$.
\end{proposition}

\proof First, we show that $\tilde\Gamma_{x,F}$  is a universal $F$--transducer.  Let $\phi^-(x,\bar y)\in L^-$ be arbitrary and let $\bar a\in F^{<\omega}$ be suitable for $\bar y$.

\begin{claim}\label{claim1}
If  $\tilde\Gamma_{x,F}(x)\wedge \phi^-(x,\bar a)$ $L$-doesn't fork over $\emptyset$ , then $\phi^-(x,\bar a)$ $L^-$-doesn't fork over $\emptyset$.
\end{claim}

\proof If $\tilde\Gamma_{x,F}(x)\wedge \phi^-(x,\bar a)$ $L$-doesn't fork over $\emptyset$, there exists $b\models\tilde\Gamma_{x,F}(x)\wedge \phi^-(x,\bar a)$ such that  $\nonforkempty{b}{\bar a}$. By Lemma \ref{lemma1}, $\minusnonforkempty{b}{\bar a}$ thus $\phi^-(x,\bar a)$ $L^-$-doesn't fork over $\emptyset$.

\begin{claim}\label{claim2}
If $\phi^-(x,\bar a)$ $L^-$-doesn't fork over $\emptyset$, then $\tilde\Gamma_x(x)\wedge \phi^-(x,\bar a)$ $L$-doesn't fork over $\emptyset$, in particular $\tilde\Gamma_{x,F}(x)\wedge \phi^-(x,\bar a)$ $L$-doesn't fork over $\emptyset$.
\end{claim}

\proof Assume  $\phi^-(x,\bar a)$ $L^-$-doesn't fork over $\emptyset$. Let $b\models \phi^-(x,\bar a)$ be such that $\minusnonforkempty{b}{\bar a}$. Let $M$ be a model of $T$. By extension in $\CC^-$, we may assume $\minusnonforkempty{b}{M\bar a}$. In particular, $tp_{L^-}(b/M\bar a)$ $L$-doesn't fork over $\emptyset$, so there exists $b^*$ such that $tp_{L^-}(b^*/M\bar a)=tp_{L^-}(b/M\bar a)$ and
$\nonforkempty{b^*}{M\bar a}$. By Remark \ref{remark1}, $b^*\models \tilde\Gamma_{x}(x)$. By the choice of $b^*$, $\phi^-(b^*, \bar a)$, thus  $\tilde\Gamma_{x}(x)\wedge \phi^-(x,\bar a)$ $L$-doesn't fork over $\emptyset$.\qed\\

\noindent It remains to show:

\begin{claim}
If $U=U(x)$ is an $\emptyset$-invariant set in $\CC$ that is an upper univesal $F$-transducer, then $U\subseteq\Gamma^*_{x,F}$. Therefore
$\tilde\Gamma_{x,F}=\Gamma^*_{x,F}$ is the greatest (with respect to inclusion) $\emptyset$-invariant set in $\CC$ that is a subset of $\CC^x$ and is a universal $F$-transducer ($\tilde\Gamma_{x,F}$ is also such greatest upper universal $F$-transducer). $\tilde\Gamma_{x,F}$ is type-definable.
\end{claim}

\proof Let $U(x)$ be as given in the claim and  assume $b\models U(x)$ and let $\nonforkempty{\bar a}{b}$. For all $\phi^-(x,\bar y)\in L^-$, if $\models\phi^-(b,\bar a)$ then
$\phi^-(x,\bar a)$ $L^-$-doesn't fork over $\emptyset$ (since $U(x)$ is an upper universal $F$-transducer). Thus $\minusnonforkempty{b}{\bar a}$, so $b\in \Gamma^*_{x,F}$.
By Lemma \ref{lemma1},  $\tilde\Gamma_{x,F}=\Gamma_{x,F}^*$ , so by Claims \ref{claim1}, \ref{claim2}, $\tilde\Gamma_{x,F}$ is the greatest $\emptyset$-invariant set in $\CC$ that is a subset of $\CC^x$ and is a universal $F$-transducer (as well as an upper universal $F$-transducer). $\tilde\Gamma_{x,F}$ is type-definable as  $\tilde\Gamma_{x,F}\equiv \bigwedge_i \Gamma_{p_i}$, where $\{p_i\}$ is the set of all complete $L$-types over $\emptyset$ of elements in $F^{<\omega}$ and $\Gamma_{p_i}$ is the partial $L$-type such that  $a\models\Gamma_{p_i}$ iff there exists $b\models p_i$ that is $
L^-$-independent from $a$ over $\emptyset$.\qed\\
$ $\\
The last statement in Proposition \ref{prop1} follows immediately by Remark \ref{trans_fk_top}. This completes the proof of Proposition \ref{prop1}.\qed\\

In order to describe the set of universal $F$-transducers, we introduce a new topology on the Stone space $S_y(T)$.

\begin{definition}\label{def. NFI-top}\em
Given a finite tuple of variables $y$, a set $U=U(y)$ is \em a basic open set in the $NFI_F$-topology \em  on $S_y(T)$ iff
there exists a type $p(x)\in S_x(T)$ with $p(x)\vdash F^{<\omega}$ and $\phi^-(x,y)\in L^-$ such that $$U=U_{p,\phi^-}=\{b\vert\ p(x)\wedge\phi^-(x,b)\ L\mbox{-doesn't fork\ over}\ \emptyset\}.$$
In case $F=\CC$, $F$ is omitted. ``$NFI$" stands for ``Non-forking instances".
\end{definition}

\begin{remark}\label {remark2}\em
As with the forking topology, we identify $\emptyset$-invariant sets with subsets of $S_y(T)$. Note that the intersection of two basic $NFI_F$-open sets is a union of basic $NFI_F$-open open sets, so the family of basic $NFI_F$-open sets forms a basis for a unique topology on $S_y(T)$.
Indeed, by extension if $b\in U_{p_0,\phi_0^-}\cap U_{p_1,\phi_1^-}$ for some $p_i,\phi^-_i$ as in Definition \ref{def. NFI-top} then
$b\in U_{q,\phi^-}$ for some $q=q(x_0,x_1)$ where $q=tp_L(a_0,a_1)$ for some independent $a_i\models p_i$ and $\phi^-=\phi^-(x_0x_1,y)=\phi^-_0(x_0,y)\wedge \phi^-_1(x_1,y)$
(clearly, $U_{q,\phi^-}\subseteq U_{p_0,\phi_0^-}\cap U_{p_1,\phi_1^-}$ and it is a basic $NFI_F$-open set). Note that since the type $p$ in Definition \ref{def. NFI-top} is a complete $L$-type, each basic $NFI_F$-open set is $L$-type-definable. Also, note that the $NFI_F$-topology will not change if we allow $p(x)$ to be a type in infinitely many variables.
\end{remark}

\begin{example}\em\label{example2}
Let $L^-=\{E\}$, $L=\{E\}\cup\{P_i\vert 0< i\leq\omega\}$ and let $T$ be the theory of an $L$-structure $M$ such that $E^M$ is an equivalence relation on its universe with infinitely many infinite $E$-classes, with exactly one class of size $i$ for every $0<i<\omega$ and such that $\{P_i^M \vert 0<i\leq\omega\}$ are pariwise disjoint and for every $0<i\leq\omega$ and $P_i^M$ is a union of exactly $i$ infinite $E$-classes. Let $T^-$ be the reduct of $T$ to $L^-$. Work in a monster model  $\CC_T$ of $T$. Now, as any $L^-$-definable set over $\emptyset$ is clearly $NFI$-open,  we conclude that each finite $E$-class is $NFI$-open. In addition, for every $0<i<\omega$, $P_i$ is a basic $NFI$-open set, while $P_\omega$ is not an $NFI$-open set (this will easily follow later, see Example-revisited \ref{ex_rev}).
\end{example}

\begin{definition}
1) A set $U\subseteq\CC$ is said to be $(L,L^-)_F$-definable over $\emptyset$ if $U=\phi^-(\CC,\bar a)$ for some $\phi^-\in L^-$ and $\bar a\in F^{<\omega}$ such that $\phi^-(x,\bar a)$
is $\emptyset$-invariant in $\CC$. If $F=\CC$ we omit $F$.\\
2) A set $U\subseteq\CC$ is said to be $(L,L^-)_F$-$\infty$-definable over $\emptyset$ if $U=p^-(\CC,\bar a)$ for some $L^-$-partial type $p^-$ over $\emptyset$ and some tuple $\bar a$ of realizations of $F$ such that $p^-(\CC,\bar a)$
is $\emptyset$-invariant in $\CC$. If $F=\CC$, we omit $F$.
\end{definition}

\begin{remark}\em\
By compactness, $U\subseteq\CC$ is $(L,L^-)_F$-$\infty$-definable over $\emptyset$ iff  $U=p^-(\CC,\bar a)$ for some $L^-$-partial type $p^-$ over $\emptyset$ and  tuple $\bar a$ of realizations of $F$ and $U$ is the solution set of
 an $L$-partial type over $\emptyset$. Likewise for $(L,L^-)_F$-definable sets over $\emptyset$.
\end{remark}

\begin{lemma}\label{lemma2}
1) If $U$ is $(L,L^-)_F$-$\infty$-definable over $\emptyset$, then $U$ is $NFI_F$-closed. If $U$ is $(L,L^-)_F$--definable over $\emptyset$, then $U$ is a basic $NFI_F$-open set.\\
2) If $T$ is stable, then $U$ is a basic $NFI_F$-open set if and only if $U$ is $(L,L^-)_F$-definable over $\emptyset$.
\end{lemma}

\proof 1) By the assumption,  there exists an $L^-$-partial type $p^-(x,\bar y)$ over $\emptyset$ and a tuple $\bar a$ (possibly infinite) of realizations of $F$ such that $U=p^-(\CC,\bar a)$ is an $\emptyset$-invariant set in $\CC$. Let $q=tp_L(\bar a)$.
Then $$p^-(\CC,\bar a)=\{b\vert\ q(\bar y)\wedge\neg\phi^-(b,\bar y)\ L\mbox{-forks\ over\ }\emptyset\ \mbox{for\ all}\ \phi^-\in p^-\}\ (*).$$ Indeed, let $R$ denote the right hand side of $(*)$. If $b\in p^-(\CC,\bar a)$ and $q(\bar y)\wedge\neg\phi^-(b,\bar y)\ L\mbox{-doesn't fork\ over\ }\emptyset$ for some $\phi^-\in p^-$ then we get contradiction to $\emptyset$-invariance of $p^-(\CC,\bar a)$ in $\CC$, so $b\in R$. If $b\not\in p^-(\CC,\bar a)$, then by $\emptyset$-invariance of $p^-(\CC,\bar a)$ in $\CC$ and extension we may assume $\nonforkempty{b}{\bar a}$. Thus $b\not\in R$. We conclude that $p^-(\CC,\bar a)$ is the intersection of complements of basic $NFI_F$-open sets. Assume now $U=\phi^-(\CC,\bar a)$ is $(L,L^-)_F$-definable over $\emptyset$. Then by $(*)$ we get immediately that $U$ is a basic $NFI_F$-open set (take $p^-(x,\bar a)=\{\neg\phi^-(x,\bar a)\}$). 2) Assume now that $T$ stable, it remains to show if $U$ is a basic $NFI_F$-open set, then it is $(L,L^-)_F$-definable over $\emptyset$. Indeed, if $U=U_{p,\phi^-}=\{b\vert\ p(x)\wedge\phi^-(x,b)\ L\mbox{-doesn't fork\ over}\ \emptyset\}$, where  $p(x)\in S_x(T)$ is such that $p(x)\vdash F^{<\omega}$ and $\phi^-(x,y)\in L^-$, then $b\in U$
iff $\phi^-(x,b)\in \bar p$ for some non-forking extension $\bar p\in S(\CC)$ of $p$. If $\bar p$ is any such extension, then there is a definition $\chi^-(y)\in L^-(\CC)$ of the $\phi^-$-type of $\bar p$ that is over $ACL^{eq}(\emptyset)$ and is a finite boolean combination of formulas of the form $\phi^-(a,y)$ for some realization $a$ of $p$ (and thus tuple of realizations of $F$) . It follows that $U=\bigvee_{i<n} \chi^-_i(\CC)$ where  $\{\chi^-_i(y)\}_{i<n}$ is the set of $\emptyset$-conjugates of $\chi^-(y)$ in $\CC$. Clearly, $U$ is $\emptyset$-invariant in $\CC$ and is an $L^-$-formula with parameters from $F$.\qed

\begin{corollary}\em
In a stable theory, a set is $(L,L^-)_F$-$\infty$-definable over $\emptyset$ iff it is a conjunction of $(L,L^-)_F$-definable sets over $\emptyset$ iff it is $NFI_F$-closed.
\end{corollary}

\proof Assume $T$ is stable. By Lemma \ref{lemma2} (1), if $U$ is $(L,L^-)_F$-$\infty$-definable over $\emptyset$ then it is $NFI_F$-closed. By Lemma \ref{lemma2} (2) an $NFI_F$-closed set is the intersection of $(L,L^-)_F$-definable sets over $\emptyset$. Finally, it is immediate that the intersection of $(L,L^-)_F$-definable sets over $\emptyset$ is $(L,L^-)_F$-$\infty$-definable over $\emptyset$.\qed

$\\$
We give now a description of the set of universal $F$-transducers via the $NFI_F$-topology.

\begin{proposition}\label{prop2}
Let $\Gamma(y)$ be an $\emptyset$-invariant set in $\CC$. Then $\Gamma(y)$ is a universal $F$-transducer iff
$\Gamma(y)$ is a dense subset of $\tilde\Gamma_{y,F}$ in the relative $NFI_F$-topology on $\tilde\Gamma_{y,F}$.
\end{proposition}

\proof By Proposition \ref{prop1}, we know that $\tilde\Gamma_{y,F}$ is a  universal  $F$-transducer and an $\emptyset$-invariant set $\Gamma=\Gamma(y)$ in $\CC$  is
an upper universal  $F$-transducer if and only if  $\Gamma\subseteq\tilde\Gamma_{y,F}$. Thus it remains to show that an $\emptyset$-invariant set $\Gamma\subseteq\tilde\Gamma_{y,F}$ in $\CC$  is a lower universal  $F$-transducer if and only if $\Gamma$ is a dense subset of $\tilde\Gamma_{y,F}$ in the relative $NFI_F$-topology on $\tilde\Gamma_{y,F}$. To show this we start with the following.

\begin{claim}\label{claim3}
For every type $p(x)\in S_x(T)$ with $p(x)\vdash F^{<\omega}$ and $\phi^-(x,y)\in L^-$, $U_{p,\phi^-}\cap \tilde\Gamma_{y,F}\neq\emptyset$ iff $\phi^-(a,y)$ $L^-$-doesn't fork over $\emptyset$ for
$a\models p$.
\end{claim}

\proof For such $p$ and $\phi^-$, $U_{p,\phi^-}\cap \tilde\Gamma_{y,F}\neq\emptyset$ iff there exists $b\models\tilde\Gamma_{y,F} $ such that $p(x)\wedge\phi^-(x,b)\ L\mbox{-doesn't fork\ over}\ \emptyset$ iff $\tilde\Gamma_{y,F}(y)\wedge\phi^-(a,y)$ $L$-doesn't fork over $\emptyset$ for $a\models p$. Since $\tilde\Gamma_{y,F}$ is a universal $F$-transducer, the latest is equivalent to $\phi^-(a,y)$ $L^-$-doesn't fork over $\emptyset$ for $a\models p$.\qed\\

\noindent Now, let $\Gamma(y)\subseteq \tilde\Gamma_{y,F}$. Then $\Gamma(y)$ is a dense subset of $\tilde\Gamma_{y,F}$ in the relative $NFI_F$-topology on $\tilde\Gamma_{y,F}$ iff
for every $p(x)\in S_x(T)$ with $p(x)\vdash F^{<\omega}$ and $\phi^-(x,y)\in L^-$ such that $U_{p,\phi^-}\cap \tilde\Gamma_{y,F}\neq\emptyset$ we have $U_{p,\phi^-}\cap\Gamma(y)\neq\emptyset$.
By Claim \ref{claim3}, the latest is equivalent to: for every $p(x)\in S_x(T)$ with $p(x)\vdash F^{<\omega}$ and $\phi^-(x,y)\in L^-$ such that $\phi^-(a,y)$ $L^-$-doesn't fork over $\emptyset$ for
$a\models p$, there exists $b\models\Gamma$ such that $p(x)\wedge\phi^-(x,b)\ L\mbox{-doesn't fork\ over}\ \emptyset$; equivalently,  for every $p(x)\in S_x(T)$ with $p(x)\vdash F^{<\omega}$ and $\phi^-(x,y)\in L^-$ such that $\phi^-(a,y)$ $L^-$-doesn't fork over $\emptyset$ for $a\models p$, the partial type $\Gamma(y)\wedge\phi^-(a,y)$ $L$-doesn't fork over $\emptyset$ for $a\models p$; namely $\Gamma(y)$ is a lower universal $F$-transducer. This completes the proof of Proposition \ref{prop2}. \qed

\begin{example-revisited}\em\label{ex_rev}
We go back to Example \ref{example2}. By Lemma \ref{lemma2}(2), it follows that a set is a basic $NFI$-open set in one variable iff it is a finite union of sets each of which is either $P_i$ for $0<i<\omega$ or
it is a finite $E$ -class. Now, if $x$ is a single variable, then easily $\tilde\Gamma_x(x)=\bigwedge_{0<i<\omega}\neg(P_i(x))$. Therefore, by Proposition \ref{prop2}, $\Gamma(x)$ is a universal transducer iff $\Gamma(x)\subseteq \bigwedge_{0<i<\omega}\neg(P_i(x))$ and $\Gamma(x)$ contains all the finite $E$-classes ($=acl_x(\emptyset)$=the set of $a\in \CC^x$ that are algebraic over $\emptyset$ in the sense of $T^-$). We conclude that there are precisely 4 universal transducers:\\ $$\tilde\Gamma_x(x)=\bigwedge_{0<i<\omega}\neg(P_i(x)),\ \tilde\Gamma_x(x)\wedge(\neg P_\omega(x)),$$  $$\tilde\Gamma_x(x)\wedge(\neg(\Lambda(x))\ \mbox{and}\ \tilde\Gamma_x(x)\wedge(\neg P_\omega(x))\wedge(\neg\Lambda(x))(=acl_x(\emptyset)),\ \mbox{where\ }$$
 $$\Lambda(x)=[\bigwedge_{0<i\leq\omega}\neg P_i(x)]\wedge (x\not\in acl_x(\emptyset)).$$
\end{example-revisited}

\begin{theorem}\label {thm}
Assume $bdd(F)=dcl^{heq}(F)$. Given variables $y$,  $\tilde\Gamma_{y,F}$ is the unique universal $F$-transducer subset of $\CC^y$ that is $(L,L^-)_F$-$\infty$-definable over $\emptyset$. Thus, if $T$ is stable, $\tilde\Gamma_{y,F}$ is the unique universal $F$-transducer subset of $\CC^y$ that is a conjunction of $(L,L^-)_F$-definable sets over $\emptyset$.
\end{theorem}

\proof First, we observe that $\tilde\Gamma_{y,F}$ is $(L,L^-)_F$-$\infty$-definable over $\emptyset$. Indeed, by Lemma \ref{lemma1}, $$\tilde\Gamma_{y,F}=\{b\in \CC^y \vert\ \minusnonforkempty{b}{bdd(F)\cap BDD(\emptyset)^-}\}.$$ For every $d\in bdd(F)\cap BDD(\emptyset)^-$, let $p_d(x,\bar f_d)=tp_{L^-}(d/\bar f_d)$, where $\bar f_d$ is a tuple
of realizations of $F$ such that $d$ is the unique solution in $\CC^{heq-}$ of $tp_{L^-}(d/\bar f_d)$ (using the assumption $bdd(F)=dcl^{heq}(F)$). Now,  $$\tilde\Gamma_{y,F}(y)=\bigwedge_{d\in D}\Lambda_d(y),\mbox{where}$$ $$\Lambda_d(y)=\exists x( p_d(x,\bar f_d)\wedge \minusnonforkempty{y}{x}),\  D=bdd(F)\cap BDD(\emptyset)^-.$$ Since each $\Lambda_d(y)$ is $L^-$-type-definable with parameters in $F$
and clearly $\tilde\Gamma_{y,F}$ is $\emptyset$-invariant in $\CC$, we get that it is $(L,L^-)_F$-$\infty$-definable over $\emptyset$. Now, let $\Gamma(y)$ be any universal $F$-transducer that is
$(L,L^-)_F$-$\infty$-definable over $\emptyset$. Then by Lemma \ref{lemma2}(1), $\Gamma(y)$ is an $NFI_F$-closed set in $S_y(T)$. By Proposition \ref{prop2}, $\Gamma(y)$ is a dense
subset of $\tilde\Gamma_{y,F}$ in the relative $NFI_F$-topology on $\tilde\Gamma_{y,F}$.  It follows that $\Gamma(y)=\tilde\Gamma_{y,F}$.\qed

\begin{remark}\em\label{neq_sorts}
All proofs in this section go through easily without the assumption that $L$ and $L^-$ have the same set of sorts; one only need to restrict the variables of $F$ and of the (upper/lower) $F$-transducers to variables  of $L^-$ and replace the universe of a model of $T$ by the universe of its restriction to $L^-$ in Remark \ref {remark1}  and Claim \ref{claim2}.
\end{remark}

\section{The lovely pair case}
Recall first the basic notions of lovely pairs. Given
$\kappa\geq |T|^{+}$, an elementary pair  $(N,M)$ of models $M\subseteq N$ of a simple theory $T$ is said to
be \em $\kappa$-lovely \em if (i) it has \em the extension property: \em for any $A\subseteq N$ of
cardinality $<\kappa$ and finitary $p(x)\in S(A)$, some nonforking extension of $p(x)$ over $A\cup
M$ is realized in $N$, and (ii) it has \em the coheir property: \em if $p$ as in (i) does not fork
over $M$ then $p(x)$ is realized in $M$. By a \em lovely pair \em (of models of $T$) we mean a
$|T|^{+}$-lovely pair.

Let $L_{P}$ be $L$ together with a new unary predicate $P$. Any elementary pair $(N,M)$ of models of $T$ ($M\subseteq N$) can
be considered as an $L_{P}$-structure by taking $M$ to be the interpretation of $P$. A basic property
from [BPV] says that any two lovely pairs of models of $T$ are elementarily equivalent, as
$L_{P}$-structures. So $T_{P}$, the common $L_{P}$-theory of lovely pairs, is complete.
$T$ has the wnfcp if every $\vert T\vert^+$-saturated model of $T_P$ is a lovely pair
(equivalently, for every $\kappa\geq\vert T\vert^+$, any $\kappa$-saturated model of $T_P$ is a $\kappa$-lovely pair).
Every theory with the wnfcp is in particular low (low theories is a subclass of simple theories).
By [BPV, Proposition 6.2], if $T$ has the wnfcp then $T_P$ is simple. Thus, this situation is a special case of our general setting in this paper, where $T_P$ is the given theory ($T$ in the general setting) and $T$ is a reduct ($T^-$ in the general setting). \\ So, in this section we assume $T$ has the wnfcp and we work in a $\lambda$-big model $\MM=(\bar M, P({\bar M}))$ of $T_P$ for some large $\lambda$   (so  $P^\MM=P(\bar M)$). $\nonforkempty{}{}$ will denote independence in $\MM$ and $\minusnonforkempty{}{}$ will denote independence in $\bar M=\MM\vert L$. Recall the following notation: for  $a\in\MM^{heq-}$, let $a^c=Cb^-(a/P(\bar M))$, where $Cb^-$ denotes the canonical base (as a hyperimaginary element) in the sense of $T$.

\begin{proposition}\label{prop3}
1) For every finite tuple of variables $x$, $\tilde\Gamma_{x}=(x=x)$, namely the greatest universal transducer in the variables $x$ is  $(x=x)$.\\
2) $\bar P(\bar x)$ and $(\neg\bar P(\bar x))\cup acl_{\bar x}(\emptyset)$ are universal transducers (where $\bar P(\bar x)$ is the conjunction $\bigwedge_i P(x_i)$, $\bar x=(x_i)_i$).\\
3) If $T$ is in addition stable (equivalently $T$ has nfcp), then the $NFI$-topology on $S_y(T_P)$ is generated by the family of $L$-definable sets over $\emptyset$. Thus an $\emptyset$-invariant set
in $\MM$ is a universal transducer iff it intersect every non-empty $L$-definable set over $\emptyset$.
\end{proposition}

We start with an observation (for part 3). Here, our notation for algebraic closure is compatible with the general setting of section 2, therefore for $A\subseteq\MM^{heq}$, $ACL^{eq}(A)$ denotes the set of imaginaries in the algebraic closure of $A$ in the sense of $\MM$ and for $A\subseteq\MM^{heq-}=\bar M^{heq}$, $acl^{eq}(A)$ denotes the set of imaginaries in the algebraic closure of $A$ in the sense of $\bar M$. We will use $tp_L(-)$ and $tp_{L_P}(-)$ for possibly hyperimaginaries in the structures $\bar M$, $\MM$ respectively.

\begin{lemma}\label{lemma3}
 $\MM^{eq-}\cap ACL^{eq}(\emptyset)=acl^{eq}(\emptyset)$ (note $\MM^{eq-}=\bar M^{eq}$).
\end{lemma}
\proof Otherwise, there exists $a\in(\MM^{eq-}\cap ACL^{eq}(\emptyset))\backslash acl^{eq}(\emptyset)$. If $a\in acl^{eq}(a^c)$, then $a\in P(\bar M)^{eq}$. Since for all
$b\in P(\bar M)^{eq}$ we have $tp_{L}(b)\equiv tp_{{L_P}}(b)$, our assumption that  $a\in(\MM^{eq-}\cap ACL^{eq}(\emptyset))$ implies $a\in acl^{eq}(\emptyset)$. So, by this  a contradiction
 we may assume $a\not\in acl^{eq}(a^c)$. By the extension property there exists a sequence
$\langle a_i \vert\ i<\omega\rangle$ of realizations of $tp_L(a/a^c)$ such that $a_0=a$ and for every $i<\omega$, $\minusnonfork{a_{i+1}}{\{a_0,...a_i\}\cup P(\bar M)}{a^c}$.
\begin{claim}
$tp_{L_P}(a_i)=tp_{L_P}(a)$ for every $i<\omega$.
\end{claim}
\proof By the construction of $\langle a_i \vert\ i<\omega\rangle$, for every $i<\omega$, $\phi^-(x,a_i)$ is realized in $P(\bar M)$ (where $x$ is a tuple of variables form the home sort of $\bar M$ and $\phi^-(x,y)\in L^{eq}$) iff  $\phi^-(x,a_i)$ $L$-doesn't fork over $P(\bar M)$ iff $\phi^-(x,a_i)$ $L$-doesn't over $a^c$ iff $\phi^-(x,a)$ $L$-doesn't fork over $a^c$ iff $\phi^-(x,a)$ $L$-doesn't fork over $P(\bar M)$ iff $\phi^-(x,a)$ is realized in $P(\bar M)$. We conclude that $Cl(tp_L(a/P(\bar M)))=Cl(tp_L(a_i/P(\bar M)))$ and thus $tp_{L_P}(a_i)=tp_{L_P}(a)$  for all $i<\omega$ (this implication is [BPV, Corollary 3.11] for real tuples but remains true for imaginary elements).\qed\\
$\\$
 Now, since $a\not\in acl^{eq}(a^c)$, we conclude that $a_{i+1}\not\in acl^{eq}(\{a_0,...a_i\})$ for all $i<\omega$ and in particular, the $a_i$-s are distinct, so $a\not\in ACL^{eq}(\emptyset)$, a contradiction. This completes the proof of Lemma \ref{lemma3} .\qed\\

\noindent\textbf{Proof of Proposition\ref{prop3}}. To prove 1),  recall the following fact (for convenience, we state it for a special case).

\begin{fact} {\em[BPV, Proposition 7.3]\em} \label{fact2}\em\\
Let $B\subseteq\MM$ and $a$ a tuple from $\MM$. Then\\ $\nonforkempty{a}{B}$ iff [$\minusnonfork{a}{B\cup P(\bar M)}{P(\bar M)}$ and $\minusnonforkempty{a^c}{B^c}$].
\end{fact}

\noindent $\Gamma^*_x=\tilde\Gamma_{x}$, so we need to show that for
every finite tuples $a,b$ from $\MM$,  $\nonforkempty{a}{b}$ implies $\minusnonforkempty{a}{b}$. By Fact \ref{fact2} it means we need to show that for
every finite tuples $a,b$ from $\MM$, if $\minusnonfork{a}{b\cup P(\bar M)}{P(\bar M)}$ and $\minusnonforkempty{a^c}{b^c}$, then $\minusnonforkempty{a}{b}$.
Indeed, as  $\minusnonfork{b}{P(\bar M)}{b^c}$, our assumption implies $\minusnonfork{b}{aP(\bar M)}{b^c}$ and in particular $\minusnonfork{b}{a}{b^c}\ (*)$. As $b^c\in dcl^{heq}( P(\bar M))$,
$\minusnonfork{a}{b^c}{a^c}$. Our assumption  $\minusnonforkempty{a^c}{b^c}$, implies $\minusnonforkempty{b^c}{aa^c}$. By (*), $\minusnonforkempty{b}{a}$.\\
\noindent We prove 2). First we show  $\bar P(\bar x)^\MM$ is a universal transducer. Assume $\phi^-(\bar x,a)$ $L$-doesn't fork over  $\emptyset$ ,where $\phi^-(\bar x,y)\in L$. By the extension property, there exists
$\bar b\in \MM$ such that $\phi^-(\bar b,a)$ and $\minusnonforkempty{\bar b}{aP(\bar M).}$ In particular, $tp_L(\bar b/aa^c)$ $L$-doesn't fork over $\emptyset$ and in particular it doesn't fork over $\bar P(\bar M)$. By the coheir property, $tp_L(\bar b/aa^c)$ is realized in $P(\bar M)$. Let $\bar b^*\in P(\bar M)$ realize it. Then $\phi^-(\bar b^*,a)$ and $\minusnonforkempty{\bar b^*}{a^c}$. By Fact \ref{fact2}, as $\bar b^*\in P(\bar M)$,  it follows that $\nonforkempty{\bar b^*}{a}$. Thus $\bar P(\bar x)\wedge \phi^-(\bar x,a)$ $L_P$-doesn't fork over $\emptyset$. By 1), we conclude that
$\bar P(\bar x)$ is a universal transducer.\\
To show that  $\Gamma(\bar x)=(\neg\bar P(\bar x))\cup acl_{\bar x}(\emptyset)$ is a universal transducer we assume $\phi^-(\bar x,a)$ $L$-doesn't fork over $\emptyset$ for $\phi^-(\bar x,y)\in L$.
If some realization of $\phi^-(\bar x,a)$ is in $acl_{\bar x}(\emptyset)$, we are done so we may assume any realization of it is not in $acl_{\bar x}(\emptyset)$. Therefore, there exists $\bar b^*\in\phi^-(\MM,a)$ such that $\minusnonforkempty{\bar b^*}{aP(\bar M)}$ and $b^*\not\in acl(aP(\bar M))$. Let $p^-=tp_L(\bar b^*/aP(\bar M))$. Let $p\in S(T_{aP(\bar M)})$ be an extension of $p^-$ that  $L_P$-doesn't fork over $\emptyset$. Let $p^*=p\vert{a}$. Then, $p^*(\bar x)\vdash (\neg\bar P(\bar x))\wedge \phi^-(\bar x,a)$, so we are done.\\
\noindent We prove 3). We need to show that for every $p(x)\in S_x(T_P)$ and $\phi^-(x,y)\in L$, the set $U_{p,\phi^-}$ is $L$-definable over $\emptyset$. We go back to the proof of Lemma \ref{lemma2} (2): Let $\chi^-(y)\in L(\MM)$ be the definition of the $\phi^-$-type of some global $L_P$-non-forking extension of $p$. Then $\chi^-(y)$ is over $ACL^{eq}(\emptyset)$. Let $c\in\MM^{eq-}$ be the canonical parameter of $\chi^-(y)$. Since $c\in ACL^{eq}(\emptyset)$, by Lemma \ref{lemma3}, $c\in acl^{eq}(\emptyset)$. As in Lemma \ref{lemma2} (2), it follows that $U_{p,\phi^-}=\bigvee_{i<n} \chi^-_i(\CC)$ where  $\{\chi^-_i(y)\}_{i<n}$ is the set of $\emptyset$-conjugates of $\chi^-(y)$ in $\MM$, but since $c\in acl^{eq}(\emptyset)$ and $acl^{eq}(\emptyset)\subseteq P(\bar M)^{eq}$, $\{\chi^-_i(y)\}_{i<n}$ is also the set of $\emptyset$-conjugates of $\chi^-(y)$ in $\MM^{eq-}=\bar M^{eq}$, so $U_{p,\phi^-}$ is $L$-definable over $\emptyset$.\qed

\begin{corollary}
Any $(L_P,L)-\infty$-definable set over $\emptyset$ containing $P(\bar x)$ must be equal to $\bar x=\bar x$.
\end{corollary}
\proof This is an immediate corollary of Theorem \ref{thm} and Proposition \ref{prop3}(1),(2).\qed


\noindent Ziv Shami, E-mail address: zivsh@ariel.ac.il.\\
Dept. of Mathematics\\
Ariel University\\
Samaria, Ariel 44873\\
Israel.\\

\end{document}